\documentclass[12pt,a4paper,twoside]{article}

\usepackage{amsmath,amsfonts,amsthm,amsopn,color,amssymb,cite,palatino}
\usepackage{graphicx}
\newcommand{\eps}{\varepsilon}
\newcommand{\C}{\mathbb{C}}
\newcommand{\R}{\mathbb{R}}
\newcommand{\Z}{\mathbb{Z}}

\newcommand{\RN}{{\mathbb{R}^N}}
\newcommand{\RT}{{\mathbb{R}^3}}

\newcommand{\RR}{{\R^+ \times \R}}
\newcommand{\de}{\partial}

\renewcommand{\le}{\leqslant}
\renewcommand{\ge}{\geqslant}
\renewcommand{\a }{\alpha }

\renewcommand{\b }{\beta }

\newcommand{\g }{\gamma }

\renewcommand{\l }{\lambda}
\newcommand{\n }{\nabla }

\renewcommand{\t}{\theta}
\renewcommand{\O}{\Omega}

\renewcommand{\S}{\Sigma}

\renewcommand{\C}{\mathbb{C}}

\newcommand{\ut}{\tilde u}
\newcommand{\vt}{\tilde v}

\newcommand{\D }{\mathcal D}
\newtheorem{theorem}{Theorem}[section]
\newtheorem{lemma}[theorem]{Lemma}

\newtheorem{definition}[theorem]{Definition}

\newtheorem{remark}[theorem]{Remark}

\renewenvironment{proof}{\noindent{\textbf{Proof\quad}}}{$\hfill\square$\vspace{0.2 cm}\\}
\newenvironment{proofmain}{\noindent{\textbf{Proof of Theorem \ref{th:main}\quad}}}{$\hfill\square$\vspace{0.2 cm}\\}
\newenvironment{proofmain2}{\noindent{\textbf{Proof of Theorem \ref{th:main2}\quad}}}{$\hfill\square$\vspace{0.2 cm}\\}
\newenvironment{proofmain3}{\noindent{\textbf{Proof of Theorem \ref{th:main3}\quad}}}{$\hfill\square$\vspace{0.2 cm}\\}

\title{Compactness results and applications to some ``zero mass'' elliptic problems}

\author{A. Azzollini \thanks{Dipartimento di Matematica, Universit\`a degli
    Studi di Bari,  Via E. Orabona 4, I-70125 Bari, Italy, e-mail: {\tt azzollini@dm.uniba.it}}
 \; \& \;A. Pomponio\thanks{Dipartimento di
    Matematica, Politecnico di Bari, Via Amendola 126/B, I-70126 Bari, Italy, e-mail: {\tt a.pomponio@poliba.it}}}

\date{}

\begin{document}

\maketitle

\begin{flushright}
{\it In ricordo di Giulio Minervini}
\end{flushright}

\

\section{Introduction and statement of the main results}

In this paper we study the elliptic problem,
\begin{equation}\label{eq0}
- \Delta v =  f'(v) \qquad \hbox{in}\;\O,
\end{equation}
in the so called ``zero mass case" that is, roughly speaking, when $f''(0)=0$.

A particular example is
\[
-\Delta v = v^{\frac{N+2}{N-2}} \qquad \hbox{in}\;\RN,
\]
with $N\ge 3$. This problem has been studied very intensely (see \cite{Au,CGS,T}) and we know the explicit expression of the positive
solutions
\[
v(x)=\frac{[N(N-2)\l^2]^{(N-2)/4}}{[\l^2+|x-x_0|^2]^{(N-2)/2}}, \qquad \hbox{with }\; \l\ge 0,\;x_0 \in \RN.
\]

If $f$ is not the critical power, we are led to require particular
growth conditions on the nonlinearity $f$. In fact, while in the
``positive mass case" (namely when $f''(0)<0$) the natural
functional setting is $H^1(\O)$ and we have suitable compact
embeddings just assuming a subcritical behavior of $f$, in the
``zero mass case" the problem is studied in $\D^{1,2}(\O)$ that is
defined as the completion of $C^\infty_0(\O)$ with respect to the
norm
\begin{equation*}
\|u\|=\left(\int_\O|\nabla u|^2\,dx\right)^\frac{1}{2}.
\end{equation*}
In order to recover analogous compactness results, we need to assume that $f$ is
supercritical near the origin and subcritical at infinity.

With these assumptions on $f$, the problem \eqref{eq0} has been
dealt with by Berestycki \& Lions~\cite{BL1,BL2,BL3}, when
$\O=\RN$, $N\ge 3$, and existence and multiplicity results have
been proved.

Recently, Benci \& Fortunato~\cite{BF} have introduced a new
functional setting, namely the Orlicz space $L^p+L^q$, which
arises very simply from the growth conditions on $f$ and seems to
be the natural framework for studying ``zero mass" problems as
shown also by Pisani in~\cite{P}.

Using this new functional setting, Benci \& Micheletti
in~\cite{BM} studied the problem \eqref{eq0}, with Dirichlet
boundary conditions, in the case of exterior domain, namely when
$\RN \setminus \O$ is contained into a ball $B_\eps.$ Under
suitable assumptions, if the ball radius $\eps$ is sufficiently
small, they are able to prove the existence of a positive
solution.

The functional setting introduced in~\cite{BF} seems to be the natural one also  for
studying the nonlinear Schr\"odinger equations with vanishing potentials, namely
\begin{equation}\label{eq:nse}
-\Delta v + V(x)v=f'(v), \qquad \hbox{in }\; \RN,
\end{equation}
with
\[
\lim_{x \to \infty} V(x)=0.
\]
Some existence results for such a problem have been found by
Benci, Grisanti \& Micheletti~\cite{BGM,BGM2} and by Ghimenti \&
Micheletti~\cite{GM}.

Even if in a different context, we need also to mention the paper of Ambrosetti, Felli
\& Malchiodi~\cite{AFM}, where problem \eqref{eq:nse} is studied when the nonlinearity
$f(v)$ is replaced by a function $f(x,v)$ of the type $K(x)v^p$, with $K$ vanishing at
infinity.

In this paper, we study problem \eqref{eq0} in two different situations.
In Section \ref{sec:prob-1}, we look for complex valued solutions of the following
problem
\begin{equation}\label{eqv}
- \Delta v =  f'(v) \qquad \hbox{in}\;\RT,
\end{equation}
assuming that $f\in C^1(\C,\R)$ satisfies the following
assumptions:
\begin{itemize}
\item[({\bf f1})] $f(0)=0$; \item[({\bf f2})] $\exists M>0$ such
that $f(M)>0;$ \item[({\bf f3})] $\forall \xi\in\C:$ $|f'(\xi)|\le
c \min(|\xi|^{p-1},|\xi|^{q-1})$; \item[({\bf f4})]
$f(e^{i\a}\rho)=f(\rho)$, for all $\xi=e^{i \a}\rho \in \C$;
\end{itemize}
where $1<p<6<q$ and $c>0$.

Observe that an example of function satisfying the previous
hypotheses can be obtained as follows. Let us consider the
function $\tilde f:\R^+ \to \R$ defined as
\[
\tilde f(t):=\left\{
\begin{array}{ll}
a t^p + b  &  \hbox{if}\;t \ge 1
\\
t^q &  \hbox{if}\;t \le 1,
\end{array}
\right.
\]
with $a,b\in \R$ chosen in order to have $\tilde f \in C^1$ and let us define
$f:\C\to \R$ as $f(\xi)=\tilde f(|\xi|)$.

Introducing the cylindrical coordinates $(r,z,\t)$, for all $n \in \Z$, we look for solutions of the type
\begin{equation}\label{eq:vn}
v^n(x,y,z)=u^n(r,z)\,e^{i n\theta}\qquad \hbox{with}\;u^n \in \R.
\end{equation}

We obtain the following existence result for problem \eqref{eqv}:
\begin{theorem}\label{th:main}
Let $f$ satisfy the hypotheses ({\bf f1}-{\bf f4}). Then
there exists a sequence $(v^n)_n$ of complex-valued solutions of problem
\eqref{eqv}, such that, for every $n\in \Z$, $v^n(x,y,z)=u^n(r,z)\,e^{i n\theta}$, with $u^n \in \R$.
\end{theorem}

Actually, an existence result in the same spirit of ours is
present in~\cite{L2}. However, in~\cite{L2} the problem is studied
using different tools and the details are omitted. Moreover in
\cite{BV} an interesting physical interpretation has been given to
the complex valued solutions of the equation \eqref{eqv} in the
positive mass case. In fact there has been shown the strict
relation between such solutions and the standing waves of the
Schr\"odinger equation with nonvanishing angular momentum.

In Section \ref{sec:prob-2}, we study
\begin{equation}\label{eqv2}
\left\{
\begin{array}{ll}
- \Delta v =  f'(v) & \hbox{in}\;\R^2\times I,
\\
v=0 & \hbox{in} \; \R^2 \times \de I,
\end{array}
\right.
\end{equation}
where $I$ is a bounded interval of $\R$ and $f\in C^1(\R,\R)$ satisfies the following assumptions:
\begin{itemize}
\item[({\bf f1'})] $f(0)=0$;
\item[({\bf f2'})] $\forall \xi\in\R:$ $f(\xi)\ge c_1  \min(|\xi|^p,|\xi|^q)$;
\item[({\bf f3'})] $\forall \xi\in\R:$ $|f'(\xi)|\le c_2 \min(|\xi|^{p-1},|\xi|^{q-1})$;
\item[({\bf f4'})] there exists $\a \ge 2$ such that $\forall \xi \in\R:\a f(\xi)\le f'(\xi)\xi$;
\end{itemize}
with $2<p<6<q$ and $c_1,$ $c_2>0.$

We will prove the following
multiplicity result:
\begin{theorem}\label{th:main2}
Let $f$ satisfy the hypotheses ({\bf f1'}-{\bf f4'}). Then there
exist infinitely many solutions with cylindrical symmetry of problem \eqref{eqv2}.
\end{theorem}

In order to approach to our problems, we use a functional
framework related to the Orlicz space $L^p+L^q$. The main
difficulty in dealing with such spaces consists in the lack of
suitable compactness results. In view of this, the key points of
this paper are two compactness theorems presented in
Section~\ref{sec:comp}. They are obtained adapting a well known
lemma of Esteban \& Lions~\cite{EL} to our situation.

The paper is organized as follows: Section~\ref{sec:L^p+L^q} is
devoted to a brief recall on the space $L^p+L^q;$ in Section~\ref{sec:comp},
we present our compactness results; in Sections~\ref{sec:prob-1} and
\ref{sec:prob-2} we solve problems \eqref{eqv} and \eqref{eqv2};
finally, in the Appendix we prove a compact embedding theorem using
similar arguments as in Section~\ref{sec:comp}.

\section{Some properties of the $L^p+L^q$
spaces}\label{sec:L^p+L^q}

In this section, we present some basic facts on the Orlicz space $L^p+L^q$. For more details, see \cite{BF,KR,P}.

Let $\O \subset \RT$. For $1<p<6<q$, denote by $(L^p(\O),\|\cdot\|_{L^p})$ and by $(L^q(\O),\|\cdot\|_{L^q})$
the usual Lebesgue spaces with their norms, and set
\begin{equation*}
L^p+L^q(\O) := \left\{v:\O\rightarrow\R\left|\:\exists(v_1,v_2)\in
L^p(\O)\times L^q(\O)\; \hbox{s.t}.\; v=v_1+v_2\right.\right\}.
\end{equation*}
The space $L^p+L^q(\O)$ is a Banach space with the norm
\begin{equation*}
\|v\|_{L^p+L^q}(\O)  :\inf\left\{\|v_1\|_{L^p}+\|v_2\|_{L^q}\mid (v_1,v_2)\in
L^p(\O) \times L^q (\O), v_1+v_2=v\right\}
\end{equation*}
and its dual is the Banach space $\big(L^{p'}(\O)\cap
L^{q'}(\O),\|\cdot\|_{L^{p'}\cap L^{q'}}\big)$, where
$p'=\frac{p}{p-1},$ $q'=\frac{q}{q-1}$ and
\begin{equation*}
\|\varphi\|_{L^{p'}\cap L^{q'}}    :\|\varphi\|_{L^{p'}}+\|\varphi\|_{L^{q'}}.
\end{equation*}

In the sequel, for all $v\in L^p +L^q(\O)$, we set
\begin{align*}
\Omega^> &:=\big\{x\in \O \mid |v(x)|>1\big\},
\\
\O^\le &:=\big\{x\in \O\mid |v(x)|\le 1\big\}.
\end{align*}

The following theorem summarizes some properties about $L^p +L^q$
spaces
\begin{theorem}
\begin{enumerate}
\item Let $v\in L^p +L^q(\O)$. Then
\begin{eqnarray}
\max\left(\|v\|_{L^q(\O^\le)}-1,
\frac{1}{1+{\rm meas}(\Omega^>)^{1/r}}\|v\|_{L^p(\Omega^>)}\right)
\phantom{\|_{L^p(\Omega_v)}).}\nonumber\\
\phantom{\max)}
\le\|v\|_{L^p+L^q}\le\max\left(\|v\|_{L^q(\O^\le)},
\|v\|_{L^p(\Omega^>)}\right)\label{controlnorm}
\end{eqnarray}
where $r=p\,q/(q-p).$
\item The space $L^p+L^q$ is continuously embedded in $L^p_{loc}$.
\item For every $r\in [p,q]:$ $L^r(\O)\hookrightarrow L^p+L^q(\O)$ continuously.
\item The embedding
\begin{equation}\label{emb}
\D^{1,2}(\O)\hookrightarrow L^p+L^q(\O)
\end{equation}
is continuous.
\end{enumerate}
\end{theorem}
\begin{proof}
\begin{enumerate}
\item See Lemma 1 in \cite{BF}.
\item See Proposition 6 of \cite{P}.
\item See Corollary 9 in \cite{P}.
\item It follows from the point 3 and the Sobolev continuous embedding
$$\D^{1,2}(\O)\hookrightarrow L^6(\O).$$
\end{enumerate}
\end{proof}

The following theorem has been proved in~\cite{P}:
\begin{theorem}\label{th:pisani}
Let $f$ be a $C^1(\C,\R)$ function (resp. $C^1(\R,\R)$) satisfying assumption ({\bf f3}) (resp. ({\bf f3'})). Then the functional
\[
v\in L^p+L^q(\O) \longmapsto \int_\O f(v) \,d x
\]
is of class $C^1$. Moreover the Nemytski operator
\[
f': v \in L^p+L^q(\O) \longmapsto f'(v) \in \left(L^p+L^q(\O)
\right)'
\]
is bounded.
\end{theorem}

Using Theorem~\ref{th:pisani} we get a very useful inequality
for the $L^p+L^q$-norm.
\begin{theorem}
For all $R>0$, there exists a positive constant $c=c(R)$ such that,
for all $v \in L^p+L^q(\O)$ with $\|v\|_{L^p+L^q} \le R$,
\begin{equation}\label{eq:dis}
\max\left(\int_{\O^>} |v|^p \,d x, \int_{\O^\le} |v|^q \,d x
\right) \le c(R)\, \|v\|_{L^p+L^q}.
\end{equation}
\end{theorem}

\begin{proof}
Let us introduce $g \in C^1(\R,\R)$ such that $g(0)=0$ and with the following growth conditions:
\begin{description}
\item[({\bf g1})] $\forall \xi\in\R:$ $g(\xi)\ge c_1  \min(|\xi|^p,|\xi|^q)$;
\item[({\bf g2})] $\forall \xi\in\R:$ $|g'(\xi)|\le c_2 \min(|\xi|^{p-1},|\xi|^{q-1})$.
\end{description}
Integrating in ({\bf g1}) we get
\[
\int_\O g(v) \, d x \ge c_1 \left( \int_{\O^>} |v|^p \,d x +
\int_{\O^\le} |v|^q \,d x  \right).
\]
By Lagrange theorem, there exists $t \in [0,1]$ such that
\[
\int_\O g'(t v)v\, d x
=\int_\O g(v) \, d x
\ge c_1 \left( \int_{\O^>} |v|^p \,d x +
\int_{\O^\le} |v|^q \,d x  \right).
\]
Then, by the boundness of $g'$ (see Theorem~\ref{th:pisani}),
there exists $M>0$ such that
\[
M \|v\|_{L^p+L^q} \ge \int_\O \left|g'(tv)v\right|\, d x \ge c_1
\left( \int_{\O^>} |v|^p \,d x + \int_{\O^\le} |v|^q \,d x \right)
\]
and hence the conclusion.
\end{proof}

\begin{remark}
Combining the inequality \eqref{controlnorm} with the estimate
\eqref{eq:dis} we deduce that the following statements are equivalent:
\begin{description}
\item{a)} $v_n\to v\hbox{ in }L^p+L^q(\O)$,
\item{b)} $\|v_n-v\|_{L^p(\O^>_n)}\to 0$ and $\|v_n-v\|_{L^q(\O^\le_n)}\to 0,$
\end{description}
where $\O^{>}_n=\{x \in \O \mid |v_n(x)-v(x)|\ge 1\}$ and $\O^\le_n$
is analogously defined.
\end{remark}

\section{Compactness results}\label{sec:comp}

In this section we present the main tools of this paper, namely a
compactness theorem for sequences with  ``a particular symmetry'' and a compact embedding of a suitable subspace of $\D^{1,2}$ into $L^p+L^q$.
The proofs of these results are both modelled on that of Theorem~1 of~\cite{EL}, which states that
a suitable subspace of $H^1$ is compactly embedded into $L^p$, for $p$ subcritical.

First of all, for every interval $I$ of $\R$, possibly unbounded, we introduce the following subspace of
$\D^{1,2}(\R^2 \times I)$:
\begin{equation*}
\D^{1,2}_{cyl}(\R^2 \times I) = \left\{u \in \D^{1,2}(\R^2 \times I) \mid u(\cdot,\cdot,z) \hbox{ is radial, for a.e. }z \in I \right\}.
\end{equation*}
Moreover we assume the following
\begin{definition}
If $u:\R^2 \times I\rightarrow\R$ is a measurable function, we
call $z$-symmetrical rearrangement of $u$ in $(x,y)$ the Schwarz
symmetrical rearrangement of the function
\begin{equation*}
u(x,y,\cdot):z\in I \mapsto u(x,y,z)\in\R.
\end{equation*}
Moreover we call $z$-symmetrical
rearrangement of $u$ the function $v$ defined as follows
\begin{equation*}
\ut:(x,y,z)\in\R^2 \times I \mapsto \ut_{x,y}(z)
\end{equation*}
where $\ut_{x,y}$ is the $z$-symmetrical rearrangement of $u$ in $(x,y).$
\end{definition}

In our first compactness result, we consider $I=\R$.
\begin{theorem}\label{th:comp}
Let $(u_j)_j$ be a bounded sequence in $\D^{1,2}_{cyl}\left( \RT \right)$ such that $u_j$ is the $z$-symmetrical
rearrangement of itself. Then $(u_j)_j$ possesses a converging subsequence in $L^p+L^q\left( \RT \right)$,
for all $1<p<6<q.$
\end{theorem}

\begin{proof}
With an abuse of notations, in the sequel for every $v\in \D^{1,2}_{cyl}\left( \RT \right)$, we denote by $v$ also the function defined in  $\RR$ as
\[
v(\sqrt{x^2+y^2},z)=v(x,y,z).
\]
Being the proof quite long and involved, we divide it into several steps, for reader's convenience.
Since $(u_j)_j$ is bounded in the $\D^{1,2}(\R^3)$ norm, there
exists $u\in \D^{1,2}_{cyl}\left( \RT \right)$ such that
\begin{eqnarray}
u_j  &  \rightharpoonup  &  u\;\hbox{ weakly in }\D^{1,2}_{cyl}\left( \RT \right) \hbox{ and in
}L^p+L^q(\R^3),\; 1< p\le 6\le q, \label{first}
\\
u_j  &  \rightarrow      &  u\;\hbox{ a.e. in }\RT, \label{second}
\\
u_j  &  \rightarrow      &  u\;\hbox{ in }L^p(K), \hbox{ for all } K\subset\subset\R^3,\;
1\le p<6.      \label{third}
\end{eqnarray}
By Lions~\cite{L},
\begin{equation}\label{decreasing}
\forall j\ge 1,\;\forall r>0,z\neq0:\quad
|u_j(x,y,z)|\le\frac{C}{r^\frac{1}{4}|z|^\frac{1}{4}},
\end{equation}
where $r=\sqrt{x^2+y^2}$.
By \eqref{decreasing}, for $R\ge 0$ large enough, $j\ge 1$ and for all $(r,|z|)\in
(R,+\infty)\times(R,+\infty)$, we have
\begin{equation}\label{eq:<1}
\begin{array}{rcl}
|u_j(x,y,z)|      &  <  &  1,\\
|u(x,y,z)|        &  <  &  1,\\
|(u_j-u)(x,y,z)|  &  <  &  1.
\end{array}
\end{equation}
Let
\begin{eqnarray*}
D_1  &  :=  &  \{(r,z)\in\R^+\times\R \mid r>R,\:|z|>R\},
\\
D_2  &  :=  &  \{(r,z)\in\R^+\times\R \mid 0\le r\le R,\:|z|\le R\},
\\
D_3  &  :=  &  \{(r,z)\in\R^+\times\R \mid 0\le r\le R,\:|z|>R\},
\\
D_4  &  :=  &  \{(r,z)\in\R^+\times\R \mid r>R,\:|z|\le R\}.
\end{eqnarray*}
Obviously $\displaystyle\bigcup_{i=1}^4 D_i=\RR$.
Moreover denote by $\chi_{D_i}$ the characteristic function of $D_i$ and observe that, since
\begin{eqnarray*}
\|u_j-u\|_{L^p+L^q}  &  =  &
\Big\| \sum_{i=1}^4(u_j-u)\chi_{D_i}\Big\|_{L^p+L^q}
\\
&  \le  &  \sum_{i=1}^4\left\| (u_j-u)\chi_{D_i}\right\|_{L^p+L^q}
=\sum_{i=1}^4\left\|u_j-u\right\|_{L^p+L^q(D_i)},
\end{eqnarray*}
then we get the conclusion if we prove that, for all $i=1,\ldots,4$,
\begin{equation*}
u_j\rightarrow u\;\hbox{ in } L^p+L^q(D_i).
\end{equation*}
\medskip

\noindent{\sc Claim 1}: \quad {\it $u_j\rightarrow u\;\hbox{ in }
L^p+L^q(D_1)$.}

Suppose for a moment that $q>8$. By \eqref{eq:<1}, for every $(x,y,z)\in D_1$, we have
$|(u_j-u)(x,y,z)|<1,$ then the inequality \eqref{controlnorm} implies
\begin{equation}\label{D1case}
\|u_j-u\|_{L^p+L^q(D_1)}\le\|u_j-u\|_{L^q(D_1)}.
\end{equation}
On the other hand, since
\[
u_j\rightarrow u \hbox{ a.e.}\quad\hbox{and }\quad
|(u_j-u)(r,z)|^q
\le\frac{C}{|r|^{\frac{q}{4}}|z|^{\frac{q}{4}}}\in L^1(D_1),
\]
by Lebesgue theorem $u_j\rightarrow u$ in $L^q(D_1)$.\\
If $6<q\le 8,$ then take $r\in \Big(q-6,4(q-6)\Big)$ and set
$\alpha=\frac{6}{6-q+r}$ and $\beta=\frac{6}{q-r}$.
Observe that $\frac{1}{\alpha}+\frac{1}{\beta}=1$ so, by Holder,
\begin{align}\label{Hold}
\int_{D_1}|u_j-u|^q\,d x\,d y\,dz  &\int_{D_1}|u_j-u|^r|u_j-u|^{q-r}\,d x\,d y\,dz\nonumber\\
&\le\Big(\int_{D_1}|u_j-u|^{\alpha r}\Big)^{\frac{1}{\alpha}}
\Big(\int_{D_1}|u_j-u|^{(\beta
q-r)}\Big)^{\frac{1}{\beta}}\nonumber\\
&\le\Big(\int_{D_1}|u_j-u|^{\frac{6 r}{6-q+r}}
\Big)^{\frac{6-q+r}{6}}\|u_j-u\|_{L^6}^{q-r}.
\end{align}
Since $(u_j)_j$ is bounded in $\D^{1,2}(\R^3)$, it is
bounded in
$L^6(\R^3).$\\
Moreover, since $q-6<r<4(q-6)$, certainly
$\frac{6r}{6-q+r}>8$, and then the last integral in
inequality \eqref{Hold} goes to zero.\\
Hence the Claim 1 is proved.\bigskip

\noindent {\sc Claim 2}:\quad {\it $u_j\rightarrow u\;\hbox{in}\;
L^p+L^q(D_2)$.}

It is enough to observe that, since $D_2$ has finite measure, $L^p+L^q(D_2)=L^p(D_2)$ (see~\cite[Remark~5]{P}) and then we get the conclusion by \eqref{third}.\bigskip

\noindent {\sc Claim 3}: \quad {\it $u_j\rightarrow u \;
\hbox{in}\; L^p+L^q(D_3)$.}

First suppose $p<4$ and consider $g\in C^1(\R,\R)$, $g(0)=0,$ such that
the following growth and strong convexity conditions hold
\begin{description}
\item{({\bf G})} $\exists k_1>0$ s.t. $\forall t\in\R:$
$|g'(t)|\le k_1\min (|t|^{p-1},|t|^{q-1})$, \item{({\bf SC})}
$\exists k_2>0$ s.t. $\forall s,t\in\R:$
$g(s)-g(t)-g'(t)(s-t)$\\
\phantom{$\exists k_2>0$ s.t. $\forall x,y\in\R:$
$g(x)-g(y)-$}$\ge k_2\min(|s-t|^p,|s-t|^q)$.
\end{description}
Since $g(0)=0$,  from ({\bf G}) and ({\bf SC}) we deduce that
\begin{equation}\label{GConf}
\exists\,k_3,k_4>0\hbox{ s.t. }\forall s\in\R:
k_3\min(|s|^p,|s|^q)\le g(s)\le k_4\min(|s|^p,|s|^q).
\end{equation}
The condition ({\bf SC}) has been introduced in~\cite{A}, where an
explicit example of function satisfying ({\bf SC}) is also given.
\\
For almost every $(x,y)\in\R^2$, we set $u^{x,y}:\R \to \R$ defined as
$u^{x,y}(z):=u(x,y,z)$. We give an analogous definition for $u^{x,y}_j$, for all $j\ge 1$.
\\
For almost every $(x,y)\in\R^2$ with $(x^2+y^2)^{1/2}\le R$, we
set
\begin{equation*}
w_j(x,y):=\int_{(-R,R)^c}g(u^{x,y}_j(z))\,dz.
\end{equation*}
We show that
\begin{equation}\label{pointconv}
w_j(x,y)\rightarrow\int_{(-R,R)^c}g(u^{x,y}(z))\,dz\quad \hbox{ for a.e. }(x,y)\in B_R.
\end{equation}
Consider
\begin{align}\label{difabsval}
\left|w_j(x,y)-\int_{(-R,R)^c}g(u^{x,y}(z))\,dz\right|
&\le\int_{(-R,R)^c}|g(u^{x,y}_j)-g(u^{x,y})|\,dz\nonumber\\
&=\int_{(-R,R)^c}|g'(\theta_j^{x,y})|\,|u^{x,y}_j-u^{x,y}|\,dz
\end{align}
where, for almost every $(x,y)\in B_R$, $\theta_j^{x,y}$
is a suitable convex combination of $u_j^{x,y}$ and
$u^{x,y}$. Since $(u_j)_j$ is bounded in $L^p+L^q$,
$g'(\theta_j^{x,y})$ is bounded in $\Big(L^p+L^q\Big)'$ (see Theorem~\ref{th:pisani}) so, by \eqref{difabsval}, to prove \eqref{pointconv}
we are reduced to show that
\begin{equation*}
u_j^{x,y}\rightarrow u^{x,y}\;\;\hbox{in} \;L^p+L^q\big((-R,R)^c\big) \;\hbox{ for a.e. }(x,y)\in B_R.
\end{equation*}
For, define
\begin{equation*}
\Omega^{x,y}_j=\left\{z \in \R \mid |z|>R,\;  |u^{x,y}_j(z)-u^{x,y}(z)|>1\right\}
\end{equation*}
so that, by \eqref{controlnorm},
\begin{multline}\label{control}
\|u^{x,y}_j-u^{x,y}\|_{L^p+L^q((-R,R)^c)}\le
\\
\max \left(\|u^{x,y}_j-u^{x,y}\|_{L^p(\Omega^{x,y}_j)},
\|u^{x,y}_j-u^{x,y}\|_{L^q((-R,R)^c \setminus \Omega^{x,y}_j)}\right).
\end{multline}
By Lebesgue theorem and by \eqref{decreasing},
\begin{equation}\label{convinLq}
\|u^{x,y}_j-u^{x,y}\|_{L^q((-R,R)^c \setminus \Omega^{x,y}_j)}\rightarrow 0\;\hbox{ for a.e. }(x,y)\in B_R.
\end{equation}
Moreover there exists $R'=R'(x,y)\in\R$ such that for all $|z|>R'$
\begin{equation*}
|u^{x,y}_j(z)-u^{x,y}(z)|\le\frac{2C}{r^{1/4}|z|^{1/4}}\le 1.
\end{equation*}
Let $\tilde R=\max(R,R')$. We have
\begin{equation*}
\|u^{x,y}_j-u^{x,y}\|^p_{L^p(\Omega^{x,y}_j)}\le\int_{(-\tilde R, -R)\cup(R,\tilde R)}|u^{x,y}_j-u^{x,y}|^p\,dz
\end{equation*}
and then
\begin{equation}\label{convinLp}
\|u^{x,y}_j-u^{x,y}\|_{L^p(\Omega^{x,y}_j)}\rightarrow 0 \;\hbox{ for a.e. }(x,y)\in B_R,
\end{equation}
because $u^{x,y}_j\to u^{x,y}$ in $L^p(\{z\in \R \mid R\le|z|\le\tilde R\})$
by \eqref{third}.
\\
By \eqref{control}, \eqref{convinLq} and \eqref{convinLp} we get
\begin{equation*}
\|u^{x,y}_j-u^{x,y}\|_{L^p+L^q((-R,R)^c)}\rightarrow 0 \;\hbox{ for a.e. }(x,y)\in B_R
\end{equation*}
and, hence, \eqref{pointconv}.
\\
We claim that the sequence $(w_j)_j$ is bounded in $W^{1,1}(B_R)$. Indeed the $L^1$-norm of $w_j$ is bounded since $(u_j)_j$ is bounded in $\D^{1,2}$ and then in $L^6$. Moreover, if we set
\begin{equation*}
\Omega_{u_j}:=\{(x,y,z) \in \RT \mid |u_j(x,y,z)|>1\},
\end{equation*}
we have
\begin{eqnarray*}
\|\nabla_{(x,y)}w_j\|_{L^1(B_R)}  &  =  &
\int_{B_R}\Big|\nabla_{(x,y)}\Big(\int_{(-R,R)^c} g(u_j(x,y,z))\,dz\Big)\Big|d x\,d y
\\
&\le & \int_{B_R}\Big(\int_{(-R,R)^c}|g'(u_j)||\nabla_{(x,y)}u_j|\,dz\Big)d x\,d y
\\
&\le & \int_{D_3}|g'(u_j)||\nabla u_j|\,d x\,d y\,dz
\\
&\le &
c_5\Big[\int_{D_3\cap\Omega_{u_j}}|u_j|^{p-1}|\nabla u_j|\,d x\,d y\,dz
\\
&     & \phantom{c}+\int_{D_3 \setminus {\Omega_{u_j}}}|u_j|^{q-1}|\nabla u_j|\,d x\,d y\,dz\Big]
\\
&\le & c_5 \Bigg[ \Big(\int_{D_3\cap\Omega_{u_j}} |u_j|^{2(p-1)}\,d x\,d y\,dz\Big)^{1/2}
\\
&     &\phantom{c_5\Big[\Big(\int_{D_3\cap\Omega_{u_j}}} \cdot
\Big(\int_{D_3\cap\Omega_{u_j}}|\nabla u_j|^2\,d x\,d y\,dz\Big)^{1/2}
\\
&     &\phantom{c}+\Big(\int_{D_3 \setminus {\Omega_{u_j}}} |u_j|^{2(q-1)}\,d x\,d y\,dz\Big)^{1/2}
\\
&     &\phantom{c_5\Big[\Big(\int_{D_3\cap\Omega_{u_j}}}
\cdot\Big(\int_{D_3 \setminus {\Omega_{u_j}}} |\nabla u_j|^2\,d x\,d y\,dz\Big)^{1/2}\Bigg]
\\
&\le & c_5\Big(\|u_j\|_{L^6}^3\|\nabla u_j\|_{L^2}+\|u_j\|_{L^6}^3\|\nabla u_j\|_{L^2}\Big)
\\
&\le & c_5\|\nabla u_j\|_{L^2}^4
\end{eqnarray*}
where we have used the fact that $2(p-1)<6<2(q-1)$ and
$\D^{1,2}(\R^3)\hookrightarrow L^6(\R^3).$ Since $W^{1,1}(B_R)$ is compactly embedded into $L^1(B_R)$,
there exists $w\in L^1(B_R)$ such that
\begin{equation}\label{eq:BR}
w_j\rightarrow w\; \hbox{ in }L^1(B_R),
\end{equation}
and, by \eqref{pointconv},
\begin{equation}\label{convL1}
w(x,y)=\int_{(-R,R)^c}g(u(x,y,z))\,dz \qquad \textrm{ a.e. in }B_R.
\end{equation}
By the definition of $w_j$ and \eqref{convL1},
\begin{align*}
\Bigg|\int_{D_3}&\Big(g(u_j(x,y,z))-g(u(x,y,z))\Big)\,d x\,d y\,dz\Bigg|
\\
&\le \int_{B_R}\left|\int_{(-R,R)^c}\Big(g(u_j(x,y,z))-g(u(x,y,z))\Big)\,dz\right| \,d x\,d y
\\
&=\|w_j -w\|_{L^1(B_R)},
\end{align*}
and so from \eqref{eq:BR} we deduce that
\begin{equation}\label{convint}
\int_{D_3}g(u_j(x,y,z))\,d x\,d y\,dz\rightarrow\int_{D_3} g(u(x,y,z))\,d x\,d y\,dz.
\end{equation}
Now observe that, by ({\bf SC}) we have
\begin{multline*}
\int_{D_3}g(u_j)-\int_{D_3}g(u)-\int_{D_3}g'(u)(u_j-u)\ge
\int_{D_{3,j}^>}|u_j-u|^p+\int_{D_{3,j}^\le}|u_j-u|^q,
\end{multline*}
where $D_{3,j}^>=\{(x,y,z)\in D_3\mid|u_j-u|>1\}$ and $D_{3,j}^\le$ is analogously defined.
\\
Moreover, by ({\bf G}) and by Proposition~29 in \cite{P},
\begin{equation*}
g'(u)\in \Big(L^p+L^q(D_3)\Big)'
\end{equation*}
and then, from \eqref{first} and \eqref{convint}, we obtain
\begin{align}
\int_{D_{3,j}^>}|u_j-u|^p & \rightarrow 0\qquad \hbox{ for}\; 1<p<4,   \label{converLp}
\\
\int_{D_{3,j}^\le}|u_j-u|^q & \rightarrow 0\qquad \hbox{for}\; q>6.  \label{converLq}
\end{align}
If $4\le p <6$, then consider $r\in \big(0,2(6-p)\big),$
$\alpha=6/(6-p+r),$ and $\beta=6/(p-r)$. Since $\frac 1\a + \frac 1\b =1$, we have
\begin{multline*}
\int_{D_{3,j}^>}|u_j-u|^p\,d x\,d y\,dz
\\
\le  \Big(\int_{D_{3,j}^>}|u_j-u|^{\alpha r}\,d x\,d y\,dz\Big)^{1/\alpha}\,
\Big(\int_{D_{3,j}^>} |u_j-u|^{(p-r)\beta}\,d x\,d y\,dz\Big)^{1/\beta}
\end{multline*}
and, since $(p-r)\beta=6$ and $\alpha r=\frac{6r}{6-p+r}<4$, from \eqref{converLp} we get
\begin{equation}\label{converLp2}
\int_{D_{3,j}^>}|u_j-u|^p\rightarrow 0\qquad \hbox{for all }4\le  p<6.
\end{equation}
From \eqref{converLp}, \eqref{converLq}, \eqref{converLp2} and using inequality \eqref{controlnorm} we have
\begin{equation*}
u_j\rightarrow u\;\hbox{ in }L^p+L^q(D_3),\qquad \textrm{ for all }1<p<6<q.
\end{equation*}

\noindent {\sc Claim 4}:\quad {\it $u_j\rightarrow u \;\hbox{ in }
L^p+L^q(D_4)$.}

The arguments are analogous to those in the previous case.\bigskip
\\
The theorem is completely proved.
\end{proof}

By the previous theorem we can easily prove the following
\begin{theorem}\label{th:comp2}
If $I\subset\R$ is bounded, then
$\D^{1,2}_{cyl}(\R^2 \times I)$ is compactly embedded in
$L^p+L^q(\R^2 \times I),$ for every $1<p<6<q$.
\end{theorem}

\begin{proof}
For the sake of simplicity, we denote $\O=\R^2 \times I$. Let $(v_j)_j\subset\D^{1,2}_{cyl}(\Omega)$ be a bounded
sequence. Up to subsequences, there exists
$v\in\D^{1,2}_{cyl}(\Omega)$ such that
\begin{align}
v_j & \rightharpoonup v \quad \hbox{ weakly in }\D^{1,2}_{cyl}(\Omega) \hbox{ and in } L^p+L^q (\O),\;1< p \le 6 \le q, \label{1convergence}
\\
v_j & \rightharpoonup v \quad\hbox{ a.e. in }\O,
\\
v_j & \rightarrow v\quad \hbox{ in } L^p(K)\quad \hbox{for all }
\,K\subset\subset\Omega,\; 1\le p<6.
\label{3convergence}
\end{align}
Let $(\hat v_j)_j$ be the sequence of the corresponding
$z$-symmetrical rearrangements of $(v_j)_j$, then we get that
there exists $w\in\D^{1,2}_{cyl}(\R^3)$ such that, up to a subsequence,
\begin{equation}\label{2convergence}
\hat v_j  \rightharpoonup w \quad\hbox{ in }\D^{1,2}_{cyl}(\R^3),
\end{equation}
and, by Theorem~\ref{th:comp},
\begin{equation}\label{4convergence}
\hat v_j  \rightarrow w \quad\hbox{ in } L^p+L^q(\R^3).
\end{equation}
We claim
\begin{description}
\item{1.} $w$ is the $z$-symmetrical rearrangement of $v$;
\item{2.} $v_j\rightarrow v$ in $L^p+L^q(\O)$.
\end{description}

1. For $R>0$, we set $B_R$ the ball in $\R^2$ of radius $R$ and centered in the origin. Let $\hat v$ be the $z$-symmetrical rearrangement of $v$. Observe that $\hat v,$ $\hat v_j\in L^p(B_R\times\R)$, indeed $v,$ $ v_j\in L^p(B_R\times I)$, for all $j \ge 1$, and
\begin{align*}
\int_{B_R\times\R}|\hat v|^p&=\int_{B_R\times I}|v|^p
\\
\int_{B_R\times\R}|\hat v_j|^p&=\int_{B_R\times
I}|v_j|^p,\quad\hbox{for all } j \ge 1.
\end{align*}
We deduce that ${\hat v}^{x,y}=\hat v(x,y,\cdot)$ and $\hat
v_j^{x,y}=\hat v_j(x,y,\cdot)$ are in $L^p(\R)$, for almost every
$(x,y)\in B_R$ and for all $j\ge 1$.
\\
Since the Schwarz symmetrization is a contraction in $L^p(\R)$
(see, for example,~\cite{AL}),
\begin{eqnarray*}
\|\hat v_j-\hat v\|^p_{L^p(B_R\times\R)} & =& \int_{B_R}\|\hat
v_j^{x,y}-\hat v^{x,y}\|^p_{L^p(\R)}\,dx\,dy
\\
& \le &
\int_{B_R}\|v_j^{x,y}-v^{x,y}\|^p_{L^p(I)}\,dx\,dy=\|v_j-v\|^p_{L^p(B_R\times
I)},
\end{eqnarray*}
therefore, by \eqref{3convergence},
\begin{equation*}
\hat v_j\rightarrow \hat v\quad\hbox{ in } L^p(B_R\times\R).
\end{equation*}
Since $R$ is arbitrary,
\begin{equation*}
\hat v_j\rightarrow \hat v\quad\hbox{ a.e. in }\R^3
\end{equation*}
so, by \eqref{4convergence}, $\hat v=w.$
\bigskip

2. Consider $g:\R\rightarrow\R$ as in the proof of Theorem~\ref{th:comp}.
Since the functional
$$u \in L^p+L^q\left( \RT \right) \longmapsto \int_{\R^3}g(u)\,dx$$
is continuous and we have proved that
\[
\hat v_j \to \hat v \quad\hbox{ in } L^p+L^q(\R^3),
\]
we get
\begin{equation}\label{6convergence}
\int_{\Omega}g(v_j)=\int_{\R^3}g(\hat v_j)\rightarrow \int_{\RT} g(\hat
v)=\int_{\Omega}g(v).
\end{equation}
Using ({\bf SC}), \eqref{1convergence}, \eqref{6convergence}
and \eqref{controlnorm}, we deduce that
\begin{equation*}
v_j\rightarrow v\quad\hbox{ in } L^p+L^q(\Omega).
\end{equation*}

\end{proof}

\section{A complex-valued solutions problem}\label{sec:prob-1}

In this section we deal with the problem
\begin{equation}\label{eqv'}
- \Delta v =  f'(v) \qquad \hbox{in}\;\RT,
\end{equation}
assuming that $f\in C^1(\C,\R)$ and satisfies the following assumptions:
\begin{itemize}
\item[({\bf f1})] $f(0)=0$; \item[({\bf f2})] $\exists M>0$ such
that $f(M)>0;$ \item[({\bf f3})] $\forall \xi\in\C:$ $|f'(\xi)|\le
c \min(|\xi|^{p-1},|\xi|^{q-1})$; \item[({\bf f4})]
$f(e^{i\a}\rho)=f(\rho)$, for all $\xi=e^{i \a}\rho \in \C$;
\end{itemize}
where $1<p<6<q$ and $c>0$.

For all $n \in \Z$, we look for solutions of the type
$v^n(x,y,z)=u^n(r,z)\,e^{i n\theta}$. Then, passing to cylindrical
coordinates, since ({\bf f4}) implies that
$f'(e^{i\a}r)=f'(r)e^{i\a}$, by some computations one can check
that, if $v^n(x,y,z)$ is solution of the problem \eqref{eqv'},
then $u^n(r,z)$ satisfies
\begin{equation}\label{eq}
- \frac{\de}{\de r}\left(r \frac{\de u^n}{\de r} \right) -
\,r\frac{\de^2 u^n}{\de z^2} +\frac{n^2}{r}\, u^n =r f'(u^n)
\qquad \hbox{in}\;\RR.
\end{equation}
Conversely, if $u^n(r,z)$ satisfies \eqref{eq}, then $v^n(x,y,z)$
is solution of \eqref{eqv'} in $\RT \setminus \R_z$, where $\R_z$
is the $z$-axis.

In the sequel we will denote with $C^{\infty}_0(\RR)$ the set of smooth functions with compact support.

Let us introduce the following Banach spaces:
\begin{itemize}
\item $L_r^s(\RR)$ the completion of $C^{\infty}_0(\RR)$ with
respect to the norm
\[
\|u\|_{L_r^s}^s:= \int_{\RR} r|u|^s d r\, d z;
\]
\item $(L^p+L^q)_r(\RR):=\{u: \RR \to \C \mid \exists (u_1,u_2)\in L^p_r \times L^q_r
\;\hbox{s.t.}\;u=u_1+u_2\}$ with the norm
\[
\|u\|_{(L^p+L^q)_r}:= \inf\{\|u_1\|_{L_r^p}+\|u_2\|_{L_r^q} \mid (u_1,u_2)\in L_r^p \times L_r^q,\;u=u_1+u_2 \};
\]
\item $E_r(\RR)$ the completion of $C^{\infty}_0(\RR)$ with respect to the norm
\[
\|u\|_r^2:\int_{\RR}
\left[r \left( \frac{\de u}{\de r}\right)^2
+r\left(\frac{\de u}{\de z}\right)^2 \right]d r\, d z;
\]
\item $E_{n,r}(\RR)$ the completion of $C^{\infty}_0(\RR)$ with respect to the norm
\[
\|u\|_{n,r}^2:\int_{\RR}
\left[r \left( \frac{\de u}{\de r}\right)^2
+r\left(\frac{\de u}{\de z}\right)^2
+ \frac{n^2}{r} u^2 \right]d r\, d z.
\]
\end{itemize}

\begin{lemma}\label{le:imm}
The following embeddings are continuous:
\[
E_{n,r}(\RR)\hookrightarrow E_{r}(\RR)\hookrightarrow (L^p+L^q)_{r}(\RR).
\]
\end{lemma}

\begin{proof}
The first one is trivial. The second one derives from the
following argument: consider the spaces
\[
(L^p+L^q)_{cyl}\left( \RT \right):=\left\{u \in L^p+L^q\left( \RT \right) \mid u(\cdot,
\cdot,z) \hbox{ is radial, for a.e. }z \in \R \right\}
\]
and
\[
{\cal F}:=\left\{u :\RT \to \R \mid u(\cdot, \cdot,z) \hbox{ is radial, for a.e. }z\in \R \right\};
\]
if we denote by $\R_z$ the $z$-axis and set $\hat
u(x,y,z)=u\big((x^2+y^2)^{\frac 12},z\big)$, by the map $u(r,z)
\to \hat u (x,y,z)$ we have
\begin{gather}
E_{r}(\RR) \simeq \overline{ C_0^\infty(\RT \setminus \R_z) \cap
{\cal F}\;}^{\;\|\cdot \|_{\D^{1,2}}} \hookrightarrow
\D^{1,2}_{cyl}(\R^2 \times \R), \label{eq:isoE}
\\
(L^p+L^q)_{r}(\RR)\simeq (L^p+L^q)_{cyl}(\R^2 \times
\R),\label{eq:isoL}
\end{gather}
so the second embedding follows immediately from \eqref{eq:isoE},
\eqref{eq:isoL} and  \eqref{emb}.
\end{proof}

For all $n \in \Z$, we will find a solution of equation \eqref{eq} looking for critical points of the functional
$J_n: E_{n,r}(\RR) \to \R$ defined as
\[
J_n(u):= \frac 12 \int_{\RR}
\left[r \left( \frac{\de u}{\de r}\right)^2
+r\left(\frac{\de u}{\de z}\right)^2
+ \frac{n^2}{r} u^2 \right]\,d r\, d z,
\]
constrained on the manifold
\[
\S_n:= \left\{ u\in E_{n,r}(\RR) \mid \int_{\RR} r f(u)\,d r\, d z=1  \right\}.
\]
By ({\bf f3}) and by Theorem~\ref{th:pisani}, the functional
\[
F_n: u\in E_{n,r}(\RR) \longmapsto  \int_{\RR} r f(u) \,d r\, d z,
\]
is well defined and continuously differentiable.

\begin{remark}
Let us observe that $\S_n$ is nonempty. Indeed, for $R>2$,
consider the functions
\begin{equation*}
\alpha_R(t):=\left\{
\begin{array}{ll}
\sqrt M(|t|-1)&\hbox{\rm if}\quad1\le |t|< 2,
\\
\sqrt M       &\hbox{\rm if}\quad2\le |t|< R,
\\
\sqrt M(R+1-|t|)&\hbox{\rm if}\quad R\le |t|< R+1,
\\
0&\hbox{\rm otherwise},
\end{array}
\right.
\end{equation*}
and
\begin{equation*}
\beta_R(t):=\left\{
\begin{array}{ll}

\sqrt M       &\hbox{\rm if}\quad1\le |t|< R,
\\
\sqrt M(R+1-|t|)&\hbox{\rm if}\quad R\le |t|< R+1,
\\
0&\hbox{\rm otherwise},
\end{array}
\right.
\end{equation*}
and set $u_R(r,z):=\alpha_R(r)\beta_R(z).$\\
Of course, $u_R \in E_{n,r}(\RR)$ and by similar arguments as in
\cite[Proof of Theorem~2]{BL1} it can be shown that, for $R$ large
enough,
\begin{equation*}
\int_{\RR} r f(u_R)\,d r\, d z>0.
\end{equation*}
Now, if $\sigma$ is a suitable rescaling parameter, the function
$$u_{R,\sigma}:(r,z)\mapsto u_R(\sigma r,\sigma z)$$
belongs to $\S_n.$
\end{remark}

A crucial step for the proof of Theorem~\ref{th:main} is the following
\begin{theorem}\label{th:min}
Let $f$ satisfy ({\bf f1}-{\bf f4}). Then, for all $n \in \Z$, there exists $u^n\in
E_{n,r}(\RR)$ such that
\[
J_n(u^n)=\min_{\S_n} J_n.
\]
\end{theorem}

\begin{proof}
Let us fix $n \in \Z$. Let $(u_j^n)_j$ be a minimizing sequence for $J_n$ constrained on $\S_n$,
namely $(u_j^n)_j$ is contained in $\S_n$ and satisfies
\[
J_n(u_j^n) \to \inf_{\S_n}J_n, \qquad \hbox{as}\;j\to \infty.
\]
Without lost of generality, we can suppose that, for all $j \ge
1$, $u_j^n$ coincides with its $z$-symmetrical rearrangement
$\ut_j^n$. Otherwise, since also $(\ut_j^n)_j$ is contained in $\S_n$ and, moreover,
\[
J_n(\ut_j^n)\le J_n(u_j^n),
\]
we should simply replace $(u_j^n)_j$ by $(\ut_j^n)_j$.
\\
Since $(u_j^n)_j$ is a bounded sequence in $E_{n,r}$, by
Lemma~\ref{le:imm} certainly there exists $u^n\in E_{n,r}$ such
that
    \begin{equation}\label{eq:inL^p+L^q_r}
        u^n_j\rightharpoonup u^n\hbox{ weakly in }E_{n,r}\hbox{ and in }
        (L^p+L^q)_r(\R_+\times\R).
    \end{equation}
On the other hand, by Lemma~\ref{le:imm} and \eqref{eq:isoE}, the
sequence $\hat u_j^n(x,y,z)=u_j^n\big((x^2+y^2)^{\frac 12},z\big)$
is bounded in $\D^{1,2}_{cyl}\left( \RT \right)$, and then, from
Theorem~\ref{th:comp}, we have that there exists $\hat u^n$ in
$(L^p+L^q)_{cyl}\left( \RT \right)$ such that, up to a
subsequence,
    \begin{equation}\label{eq:inL^p+L^qcyl}
        \hat u_j^n\to\hat u^n\hbox{ in }(L^p+L^q)_{cyl}\left( \RT
        \right).
    \end{equation}
By continuity,
\begin{equation}\label{eq:intconv}
\int_\RT f(\hat u_j^n)\to\int_\RT f(\hat u^n), \qquad \hbox{as}\;j \to \infty.
\end{equation}\label{intconv}
Moreover, comparing \eqref{eq:inL^p+L^q_r} and
\eqref{eq:inL^p+L^qcyl}, by \eqref{eq:isoL} we deduce that $\hat
u^n(x,y,z)=u^n\big((x^2+y^2)^{\frac 12},z\big),$ so, passing to
cylindrical coordinates in \eqref{eq:intconv}, we have
\[
\int_{\RR}r f(u^n_j) \,d r\,dz \to \int_{\RR}r f(u^n) \,d r\,dz,
\qquad \hbox{as}\;j \to \infty,
\]
and then $u^n \in \S_n$.\\
Finally, the weak lower semicontinuity of the $E_{n,r}$-norm and
\eqref{eq:inL^p+L^q_r} imply
\[
J_n(u^n)=\min_{\S_n}J_n,
\]
and the theorem is proved.
\end{proof}

Now we show how Theorem~\ref{th:main} follows immediately from
Theorem~\ref{th:min}.

\begin{proofmain}
Let us fix $n\in \Z$.
Let $u^n \in E_{n,r}(\RR)$ be a minimizer of $J_n$ constrained on $\S_n$,
whose existence is guaranteed by Theorem~\ref{th:min}.
Since $u^n$ is a critical point of $J_n$ constrained on $\S_n$, there exists
a Lagrange multiplier $\l^n$ such that $(\l^n,u^n)$ satisfies
\begin{equation*}
- \frac{\de}{\de r}\left(r \frac{\de u^n}{\de r} \right)
- \,r\frac{\de^2 u^n}{\de z^2}
+\frac{n^2}{r}\, u^n
=\l^n \, r f'(u^n)
\qquad \hbox{in}\;\RR,
\end{equation*}
and so, as already observed
\begin{equation*}
- \Delta v^n = \l^n \,f'(v^n) \qquad \hbox{in}\;\RT \setminus \R_z,
\end{equation*}
for $v^n(x,y,z)=u^n(r,z) e^{in\t}$. Following the idea
of~\cite{BV} (see Theorem~3, therein), we argue that $(\l^n,v^n)$
is, in fact, a solution of
\begin{equation*} 
- \Delta v^n = \l^n \,f'(v^n) \qquad \hbox{in}\;\RT.
\end{equation*}
Let us observe that $\l^n>0$. Indeed, by ({\bf f4}), we can find $w^n\in
\D^{1,2}\left( \RT \right)$, of the type $w^n(x,y,z)=\g^n(r,z)e^{in\t}$, with $\g^n(r,z)\in \R$, such
that
\[
\int_\RT f'(v^n) {\bar w^n} d x\, d y\, d z= \int_\RT f'(u^n) {\g}^n d x\, d y\, d z>0.
\]
Now we can repeat the arguments of \cite[Theorem~2]{BL1} to prove that $\l^n>0$.
\\
Finally, it is easy to see that
\[
v^n_{\l^n}(x):=v^n\left(\frac{1}{\sqrt{\l^n}}\, x\right)
\]
is a solution of \eqref{eqv'}.
\end{proofmain}

\section{A ``zero mass" problem in $\R^2 \times
I$}\label{sec:prob-2}

In this section we consider the problem
\begin{equation}\label{eqv2'}
\left\{
\begin{array}{ll}
- \Delta v =  f'(v) & \hbox{in}\;\R^2\times I,
\\
v=0 & \hbox{in} \; \R^2 \times \de I,
\end{array}
\right.
\end{equation}
where $I$ is a bounded interval of $\R$ and $f\in C^1(\R,\R)$
satisfies the following assumptions:
\begin{itemize}
\item[({\bf f1'})] $f(0)=0$;
\item[({\bf f2'})] $\forall \xi\in\R:$ $f(\xi)\ge c_1  \min(|\xi|^p,|\xi|^q)$;
\item[({\bf f3'})] $\forall \xi\in\R:$ $|f'(\xi)|\le c_2
\min(|\xi|^{p-1},|\xi|^{q-1})$;
\item[({\bf f4'})] there exists $\a \ge 2$ such that $\forall \xi \in\R:\a f(\xi)\le f'(\xi)\xi$;
\end{itemize}
with $2<p<6<q$ and $c_1,\,c_2>0$.
\\
Set $\Omega:=\R^2\times I.$
By $({\bf f3'})$ and Theorem~\ref{th:pisani}, the functional
\begin{equation*}
J(v):=\frac{1}{2}\int_{\O}|\nabla v|^2\,dx-\int_{\O}f(v)\,dx,\quad v\in\D^{1,2}(\Omega),
\end{equation*}
is $C^1$ and its critical points are weak solutions of
\eqref{eqv2'}.

In particular, we are interested in finding solutions with
cylindrical symmetry. Since $\D^{1,2}_{cyl}(\Omega)$ is a natural
constraint, that is every critical point of the functional $J$
constrained on $\D^{1,2}_{cyl}(\Omega)$ is a critical point of the
nonconstrained functional, we will look for critical points of
$J_{|\D^{1,2}_{cyl}(\Omega)}.$ To simplify the notations, from now
on we set $\widehat J=J_{|\D^{1,2}_{cyl}(\Omega)}.$

Now we are ready to prove Theorem~\ref{th:main2}.

\begin{proofmain2}
Since $\widehat J$ is even and of class $C^1$ , we may apply a
very well known symmetrical version of the mountain pass theorem
(see~\cite{AR} or \cite{BBF}). We just have to verify the following three conditions:
\begin{description}
\item{1.} $J$ satisfies the Palais-Smale condition, i.e.
any sequence $(v_j)_j$ in $\D^{1,2}_{cyl}(\O)$ such that
\begin{equation}\label{eq:PS}
\big(\widehat J(v_j)\big)_j \hbox{ is bounded,}
\qquad
\widehat J'(v_j)\rightarrow 0,
\end{equation}
admits a convergent subsequence;
\item{2.} there exist $\rho>0$ and $C>0$ such that
\begin{equation*}
\widehat J(u)>C,\quad \hbox{for all }u\in S_\rho,
\end{equation*}
where $S_\rho:=\big\{u\in\D_{cyl}^{1,2}(\O)\mid \|u\|=\rho\big\};$
\item{3.} for all $V\subset \D^{1,2}_{cyl}(\O)$ such that $\dim V<+\infty$, we have
$\displaystyle\lim_{\substack{ u\rightarrow +\infty \\ u\in V}}\widehat J(u)=-\infty.$
\end{description}
For the proof of the $2^{nd}$ and $3^{rd}$ conditions, we refer to
\cite[Propositions~33 and 34]{P}.
\\
As regards the Palais-Smale condition, we first observe
that, by standard arguments, the hypotheses \eqref{eq:PS}
imply that the sequence $(v_j)_j$ is bounded in
$\D^{1,2}_{cyl}(\O).$ So there exists
$v\in\D^{1,2}_{cyl}(\O)$ such that, up to a subsequence,
\begin{equation*}
v_j\rightharpoonup v\quad\hbox{weakly in }
\D^{1,2}_{cyl}(\O),
\end{equation*}
and, by Theorem~\ref{th:comp2},
\begin{equation}\label{eq:strconv}
v_j\to v\quad\hbox{in }L^p+L^q(\O).
\end{equation}
Now, since
$$\widehat J'(v_j)=-\Delta v_j-f'(v_j),$$
from the second of \eqref{eq:PS}, we deduce that there exists an infinitesimal sequence $(\varepsilon_j)_j$ such that
$$-\Delta v_j=f'(v_j)+\varepsilon_j\quad\hbox{in }\big(\D^{1,2}_{cyl}(\O) \big)'.$$
Thus, inverting the Laplacian and using \eqref{eq:strconv}
and the continuity of the Nemytski operator associated with $f'$, we have
$$v_j=(-\Delta)^{-1}(f'(v_j)+\varepsilon_j)\to(-\Delta)^{-1}(f'(v)),$$
and we are done.
\end{proofmain2}

\appendix
\section{Appendix}

This section is entirely devoted to the proof of the following
compact embedding theorem
\begin{theorem}\label{th:main3}
Let $N=\sum_{i=1}^m N_i$, with $N \ge 3$, $m\ge 1$ and $N_i\ge 2$ for all $1\le i\le m$.
Then the space
\begin{align*}
\D_s^{1,2}(\RN):=\Big\{u\in\D^{1,2}(\RN)\mid & \;\forall i\in\{1,\ldots,m\},
\forall x_i^0 \in \R^{N_i},
\\
& u(x^0_1,\ldots,x^0_{i-1},\,\cdot\, ,x^0_{i+1},\ldots,x^0_m) \hbox{\rm{ is radial}} \Big\}
\end{align*}
is compactly embedded in $L^p+L^q(\R^N),$ for all $1<p<2^*<q$, where $2^*=2N/(N-2)$.
\end{theorem}

First of all, let us observe that, when $m=1$, Theorem~\ref{th:main3} has been already proved in~\cite{BF}. In this paper, we will deal with the case $m \ge 2$. This theorem is used in~\cite{A} to find solutions for the semilinear Maxwell equations in even dimension.

\begin{remark}
Let $H^1_s(\RN) \subset H^1(\RN)$ be defined in an analogous way. In~\cite{L3}, Lions proved that $H^1_s(\RN)$ is compactly embedded into $L^p(\RN)$ for $2<p<2^*$.
\end{remark}

To prove Theorem~\ref{th:main3} first we need to introduce two preliminary lemmas.
\begin{lemma}
Let $u \in \D_s^{1,2}(\RN)$ be such that $u$ is decreasing with respect to $|x_i|$, for $i \ge 2$. Then
\begin{multline}\label{eq:stima}
0 \le u \le C_N \left(
\|u\|_{L^{2^*}}^{N/(2N-2)}\,
\|\n_{x_1} u\|_{L^2}^{(N-2)/(2N-2)}\right) \cdot
\\
\cdot |x_1|^{-(N_1 -1)(N-2)/(2N-2)}
\prod_{i=2}^m |x_i|^{-N_i (N-2)/(2N-2)},
\end{multline}
where $C_N>0$ depends only on $N$.
\end{lemma}

\begin{proof}
The proof is simply a combination of \cite[Proposition~2.1]{L} and \cite[Lemma~3.1]{L3}.
\end{proof}

\begin{lemma}\label{le:compemb}
Let $m\ge 1$, $n \ge 2$ and $N=m+n$. Let $(u_j)_j$ be a bounded sequence in $\D^{1,2}(\R^m \times \R^n)$, such that, for all $j \ge 1$ and for almost every $x\in \R^m$, the function $u_j(x,\cdot)$ is radial in $\R^n$. Then there exists $u\in \D^{1,2}(\R^m \times \R^n)$ such that for all $\O\subset\R^m$ bounded and with Lipschitz
boundary $(u_j)_j$ converges, up to subsequences, to $u$ in $L^p+L^q(\O \times \R^n)$, for $1<p<2^*<q$.
\end{lemma}

\begin{proof}
Since $(u_j)_j$ is a bounded sequence in $\D^{1,2}(\R^m \times \R^n),$ there exists $u\in\D^{1,2}(\R^m \times \R^n)$ such that
\begin{align}
u_j\rightharpoonup &u\;\hbox{ weakly in }\D^{1,2}(\R^m \times \R^n),\nonumber
\\
u_j\to &u\;\hbox{ a.e. in }  \R^m \times \R^n.\label{eq:q.o.}
\end{align}
Following an idea of \cite{ABDF}, we set
\begin{equation}
w_j=\left\{
\begin{array}{ll}
|u_j-u| &\hbox{if }|u_j-u|\ge\varepsilon,
\\
\eps^{-2}|u_j-u|^3&\hbox{if }|u_j-u|\le\varepsilon,
\end{array}
\right.
\end{equation}
with $0<\eps<1$. Observe that $w_j$ is well defined and almost everywhere differentiable.
By some computations we get
\begin{align*}
|w_j|^2&\le \eps^{-4}|u_j-u|^6,
\\
|\nabla w_j|^2 &\le 9|\nabla u_j- \n u|^2,
\end{align*}
so the sequence $(w_j)_j$ is bounded in $H^1(\R^m\times\R^n).$
\\
Now, fix $\O\subset\R^m$ bounded, with Lipschitz boundary and, with an abuse of notations,
relabel by $w_j,$ $u_j$ and $u$ the restrictions of the same
functions to $\O\times\R^n.$ By \cite[Lemma~3.2]{L3}, there exists
$w\in H^1(\O\times\R^n)$ such that
\begin{align*}
w_j&\to w\quad\hbox{in }L^r(\O\times\R^n),\hbox{ for }2<r<2^*,
\\
w_j&\to w\quad\hbox{a.e. in }\O\times\R^n.
\end{align*}
On the other hand, by \eqref{eq:q.o.}, we infer that $w=0$. Let us consider for a moment $p>2$. As in \cite{ABDF}, we get
\begin{align*}
\int_{\{|u_j- u|\ge 1\}} &|u_j-u|^p \, dx
+\int_{\{|u_j- u|\le 1\}} |u_j-u|^q \, dx
\\
&\le 2\int_{\{|u_j- u|\ge \eps \}} |u_j-u|^p \, dx
+\int_{\{|u_j- u|\le \eps\}} |u_j-u|^q \, dx
\\
&\le 2\int_{\{|u_j-u|\ge \eps \}} |w_j|^p \, dx
+\eps^{q-2^*}\int_{\{|u_j- u|\le \eps\}} |u_j-u|^{2^*} \, dx
\\
& \le 2 \|w_j\|_{L^{p}(\O\times\R^m)}^p + \eps^{q-2^*} \|u_j-u\|_{L^{2^*}(\O\times\R^m)}^{2^*}.
\end{align*}
If, instead, $1<p\le 2$, taking $2<r<2^*$, we have:
\begin{align*}
\int_{\{|u_j- u|\ge 1\}} &|u_j-u|^p \, dx
+\int_{\{|u_j- u|\le 1\}} |u_j-u|^q \, dx
\\
&\le \int_{\{|u_j- u|\ge 1\}} |u_j-u|^r \, dx
+\int_{\{|u_j- u|\le 1\}} |u_j-u|^q \, dx
\\
& \le 2 \|w_j\|_{L^{r}(\O\times\R^m)}^r + \eps^{q-2^*} \|u_j-u\|_{L^{2^*}(\O\times\R^m)}^{2^*}.
\end{align*}
Therefore, in any case, we can conclude that for all $1<p<2^*$, we have
\begin{multline*}
\int_{\{|u_j- u|\ge 1\}} |u_j-u|^p \, dx
+\int_{\{|u_j- u|\le 1\}} |u_j-u|^q \, dx
\\
\le 2 \|w_j\|_{L^{r}(\O\times\R^m)}^r + \eps^{q-2^*} \|u_j-u\|_{L^{2^*}(\O\times\R^m)}^{2^*}.
\end{multline*}
with a suitable $2<r<2^*$.
\\
Since $w_j\to0$ in $L^r(\O\times\R^m)$, for all $2<r<2^*$, and $(u_j)_j$ is bounded
in $L^{2^*}(\O\times\R^m)$, by the arbitrariness of $\eps$, we
infer that
\begin{equation*}
\int_{\{|u_j- u|\ge 1\}} |u_j-u|^p \, dx+
\int_{\{| u_j-u|\le 1\}} | u_j-u|^q \, dx \to0.
\end{equation*}
Thus the conclusion follows from \eqref{controlnorm}.
\end{proof}

Now we pass to prove Theorem~\ref{th:main3}.

\begin{proofmain3}
Let $(u_j)_j$ be a bounded sequence in $\D^{1,2}_s(\RN)$. Up to a subsequence, we have
\begin{align*}
u_j &\rightharpoonup u \quad\hbox{weakly in }\D^{1,2}_s(\RN),\\
u_j &\to u\quad\hbox{in } L^p(K),\quad\forall K\subset\subset\RN,\quad 1\le p< 2^*,\\
u_j &\to u\quad\hbox{a.e. in }\RN.
\end{align*}
Let $v_j=S_2(S_3(\ldots(S_m(u_j))\ldots)),$ where $S_i$ is the symmetrizing operator with respect to $x_i \in \R^{N_i}$, with $i=2,\ldots,m$.

\medskip

\noindent {\sc Step 1}: \quad {\it There exists $v \in
\D^{1,2}_s(\RN)$ such that, up to a subsequence, $v_j \to v$ in
$L^p+L^q(\RN)$.}

We will just sketch the proof, since it can modelled on that of Theorem~\ref{th:comp}.
\\
First of all, since $(v_j)_j$ is bounded in $\D^{1,2}_s(\RN)$, there exists $v \in \D^{1,2}_s(\RN)$ such that
\begin{align*}
v_j &\rightharpoonup v \quad\hbox{weakly in }\D^{1,2}_s(\RN),\\
v_j &\to v\quad\hbox{in } L^p(K),\quad\forall K\subset\subset\RN,\quad 1\le p< 2^*,\\
v_j &\to v\quad\hbox{a.e. in }\RN.
\end{align*}
Let us observe that by \eqref{eq:stima}, there exists  $R>0$ such that if $|x_i|>R$, for all $i=1,\ldots, m$, then
\begin{eqnarray*}
\big|v_j(x_1,\ldots,x_m) -v(x_1,\ldots,x_m)\big|&<&1,
\\
\big|v_j(x_1,\ldots,x_m)\big|&<&1,
\\
\big|v(x_1,\ldots,x_m)\big|&<&1.
\end{eqnarray*}
Arguing as in the proof of Theorem~\ref{th:comp}, we need only to check that
\[
v_j \to v\qquad \hbox{in }L^p+L^q (D_{I_1,I_2}),
\]
where
\[
D_{I_1,I_2}:=\big\{x=(x_1,\ldots,x_m)\in \RN \mid |x_i|\ge R,\hbox{ if }i\in I_1;\;
|x_i|\le R,\hbox{ if }i\in I_2\big\},
\]
for every $I_1$ and $I_2$ such that $I_1\cap I_2 =\emptyset$ and $I_1\cup I_2 =\{1,\ldots,m\}$.
\\
Let us define, in particular, the following sets:
\begin{align*}
D':=&\big\{x\in \RN \mid |x_i|\ge R,\hbox{ for all }i=1,\ldots,m\big\},
\\
D'':=&\big\{x\in \RN \mid |x_i|\le R,\hbox{ for all }i=1,\ldots,m\big\},
\\
D''':=&\big\{x\in \RN \mid |x_i|\le R,\hbox{ if }i=1,\ldots,k,\;|x_i|\ge R,\hbox{ if } i=k+1,\ldots,m\big\}.
\end{align*}
Without lost of generality, we need to prove the
$L^p+L^q$-convergence of $(v_j)_j$ only in these three particular
domains.

Arguing as in the Claim 1 of the proof of Theorem~\ref{th:comp}, we can prove that $v_j \to v$ in $L^p+L^q(D')$. Indeed, if $q>\frac{N_1(2N-2)}{(N_1-1)(N-2)}$, the convergence follows immediately by \eqref{eq:stima} and by Lebesgue theorem. If, instead, $2^*<q\le \frac{N_1(2N-2)}{(N_1-1)(N-2)}$, then take
\[
r \in \left(q -2^*, \frac{N_1(N-1)(q -2^*)}{N-N_1}\right),
\]
and set $\alpha=\frac{2^*}{2^*-q+r}$ and $\beta=\frac{2^*}{q-r}$.
Observe that $\frac{1}{\alpha}+\frac{1}{\beta}=1$ so, by Holder,
\begin{align}
\int_{D'}|v_j-v|^q\,d x &\int_{D'}|v_j-v|^r|v_j-v|^{q-r}\,d x\nonumber
\\
&\le\Big(\int_{D'}|v_j-v|^{\alpha r}\Big)^{\frac{1}{\alpha}}
\Big(\int_{D'}|v_j-v|^{\beta
q-r)}\Big)^{\frac{1}{\beta}}\nonumber
\\
&\le\Big(\int_{D'}|v_j-v|^{\frac{2^* r}{2^*-q+r}} \Big)^{\frac{2^*-q+r}{2^*}}\|v_j-v\|_{L^{2^*}}^{q-r}. \label{eq:H}
\end{align}
Since $(v_j)_j$ is bounded in $L^{2^*}(\RN)$ and $\frac{2^*-q+r}{2^*}>\frac{N_1(2N-2)}{(N_1-1)(N-2)}$, the last integral in \eqref{eq:H} goes to zero.

Being $D''$ a finite measure set, we deduce easily that $v_j \to
v$ in $L^p+L^q(D'')$.

Let us now consider the case of $D'''$. If $p<\frac{2N-2}{N-2}$, let $g$ be a $C^1$-function with $g(0)=0$ and  satisfying the following conditions:
\begin{description}
\item{({\bf G})} $\exists c_1>0$ s.t. $\forall t\in\R:$
$|g'(t)|\le c_1\min (|t|^{p-1},|t|^{q-1})$, \item{({\bf SC})}
$\exists c_2>0$ s.t. $\forall s,t\in\R:$
$g(s)-g(t)-g'(t)(s-t)$\\
\phantom{$\exists c_2>0$ s.t. $\forall x,y\in\R:$
$g(x)-g(y)-$}$\ge
c_2\min(|s-t|^p,|s-t|^q)$.
\end{description}
Moreover, for all $\bar x=(x_1,\ldots,x_k)$, we set
\[
w_j(\bar x)=\int_{|x_{k+1}|\ge R}\cdots \int_{|x_{m}|\ge R}g(v_j(\bar x,x_{k+1},\ldots,x_m))\,d x_{k+1} \cdots d x_m.
\]
Following the scheme of Claim 3 in the proof of Theorem~\ref{th:comp}, we show that $v_j \to v$ in $L^p+L^q(D''')$.
Finally, if $1<p\le \frac{2N-2}{N-2}$, we get the conclusion again
by Holder inequality.

Therefore Step 1 is completely proved.

\medskip

\noindent {\sc Step 2}: \quad {\it $u_j\to u$ in $L^p+L^q(\RN)$.}

Define $\vt_j=S_3(S_4(\ldots(S_m(u_j))\ldots))$. Observe that $v_j=S_2(\vt_j)$. Consider $R>0$ and
\[
Q_R=\!\{x=(x_1,\ldots,x_m)\in \RN \!\mid |x_1|\le R,\,x_2\in \R^{N_2},\,|x_3|\le R,\ldots,|x_m|\le R\}.
\]
By Lemma~\ref{le:compemb}, we have that there exists $\vt \in \D^{1,2}(\RN)$ such that, up to subsequences,
\begin{align}
\vt_j &\rightharpoonup \vt \quad\hbox{weakly in }\D^{1,2}_s(\RN) \hbox{ and in } L^p+L^q(\RN), \label{eq:2''}
\\
\vt_j &\to \vt \quad\hbox{in } L^p+L^q(Q_R).    \label{eq:3}
\end{align}
Let us show that $S_2(\vt)=v$. Let $g\in C^1(\R,\R)$ be even and satisfying the growth conditions ({\bf G}) and ({\bf SC}). By \cite[Corollary~2.3]{AL}, we have
\begin{align*}
\int_{Q_R} g\big(S_2(\vt_j) -S_2(\vt)\big) &
= \int_{Q_R} g\Big(\big|S_2(\vt_j) -S_2(\vt)\big|\Big)
\\
&=\int_{Q_R} g\Big(\big|S_2(|\vt_j|) -S_2(|\vt|)\big|\Big)
\\
&\le \int_{Q_R} g\Big(\big||\vt_j| -|\vt|\big|\Big)
\\
& \le \int_{Q_R} g\big(|\vt_j -\vt|\big)
= \int_{Q_R} g(\vt_j -\vt).
\end{align*}
By \eqref{eq:3}, ({\bf G}) and Theorem \ref{th:pisani}  we deduce
that
\[
\int_{Q_R} g(\vt_j -\vt) \to 0.
\]
Therefore
\[
\int_{Q_R} g\big(S_2(\vt_j) -S_2(\vt)\big) \to 0
\]
and then, by \eqref{eq:dis},
\[
v_j= S_2(\vt_j) \to S_2(\vt)\quad\hbox{in } L^p+L^q(Q_R).
\]
By Step 1 and by the arbitrariness of $R>0$, we infer that $S_2(\vt)=v$.
\\
Hence
\begin{equation*}
\int_\RN g(\vt_j)
=\int_\RN g(v_j) \to
\int_\RN g(v)
=\int_\RN g(\vt),
\end{equation*}
so, by \eqref{eq:2''} and ({\bf SC}), we conclude that
\[
\vt_j \to \vt \quad\hbox{in } L^p+L^q(\RN).
\]
Iterating this argument we show that
$v=S_2(S_3(\ldots(S_m(u))\ldots))$ and hence the conclusion.

\end{proofmain3}

\end{document}